\documentstyle[fleqn,12pt]{article}
\begin{document}\pagenumbering{arabic}\setcounter{page}{1}
\pagestyle{plain}
\baselineskip=16pt

\thispagestyle{empty}
\rightline{MSUMB 98-03, December 1998} 
\vspace{1.4cm}

\begin{center}
{\Large\bf Differential Geometry of the $q$-plane}
\end{center}

\vspace{1cm}
\begin{center} Sultan A. \c Celik \footnote{E-mail: celik@yildiz.edu.tr}
and Salih \c Celik \footnote{E-mail: sacelik@yildiz.edu.tr}\\
Yildiz Technical University, Department of Mathematics, \\
80270 Sisli, Istanbul, TURKEY. \end{center}

\vspace{2cm}
{\bf Abstract}

Hopf algebra structure on the differential algebra of the extended $q$-plane 
is defined. An algebra of forms which is obtained from the 
generators of the extended $q$-plane is introduced and its Hopf algebra 
structure is given. 

\vfill\eject
\noindent
{\bf 1. Introduction}

\noindent
Differential geometry of Lie groups plays an important role in the 
mathematical modelling of physics theories. A class of noncommutative Hopf 
algebra has been found in the discussions of the integrable systems. These 
Hopf algebras are $q$-deformed function algebras of classical groups and 
this structure is called quantum group.$^1$ The quantum group can also be 
regarded as a generalization of the notion of a group.$^2$ Thus it is also 
attracting to generalize the corresponding notions of differential geometry. 
Mathematical aspects of such a generalization are promising. 

Noncommutative geometry$^3$ is one of the most attractive mathematical 
concepts in physics and has started to play an important role 
in different fields of mathematical physics for the last few years. 
The basic structure giving a direction to the noncommutative geometry 
is a differential calculus$^4$ on an associative algebra. Quantum plane 
is simple example of quantum space, and to investigate its 
differential geometry may be interesting which is considered in the 
present work. 

\noindent
{\bf 2. Review of Hopf algebra ${\cal A}$} 

\noindent
{\bf 2.1 The algebra of polynomials on the $q$-plane}

\noindent
The quantum plane or the Manin plane $^5$ is defined as an associative 
algebra generated by two noncommuting coordinates $x$ and $y$ with the 
single quadratic relation 
$$ x y - q y x = 0, \qquad q \in {\mbox C} - \{0\}. \eqno(1) $$
This associative algebra over the complex number, C, is known 
as the algebra of polynomials over quantum plane and is often denoted by 
${\cal C}_q[x,y]$. We shall just call it ${\cal C}_q$. In the limit 
$q \longrightarrow 1$, this algebra is commutative and can be considered 
as the algebra of polynomials ${\cal C}[x,y]$ over the usual plane, where 
$x$ and $y$ are the two coordinate functions. We denote the unital extension 
of ${\cal C}_q$ by ${\cal A}$, i.e. it is obtained by adding a unit element. 

\noindent
{\bf 2.2 Hopf algebra structure on ${\cal A}$}

\noindent
Hopf algebra structure on the extended quantum plane was introduced 
in ref. 6. The definitions of a coproduct and a counit on the algebra 
${\cal A}$ as follows: 

{\bf (1)} The C-algebra homomorphism (coproduct) 
$\Delta: {\cal A} \longrightarrow {\cal A} \otimes {\cal A}$ 
is given by 
$$\Delta(x) = x \otimes x, $$
$$\Delta(y) = y \otimes 1 + x \otimes y, \eqno(2)$$
which is coassociative: 
$$(\Delta \otimes \mbox{id}) \circ \Delta = 
  (\mbox{id} \otimes \Delta) \circ \Delta \eqno(3)$$
where id denotes the identity map on ${\cal A}$. 

{\bf (2)} The C-algebra homomorphism (counit) 
$\epsilon: {\cal A} \longrightarrow {\mbox C}$ is defined by 
$$\epsilon(x) = 1, \qquad \epsilon(y) = 0. \eqno(4)$$
The counit $\epsilon$ has the property 
$$\mu \circ (\epsilon \otimes \mbox{id}) \circ \Delta 
  = \mu' \circ (\mbox{id} \otimes \epsilon) \circ \Delta \eqno(5)$$
where $\mu : {\mbox C} \otimes {\cal A} \longrightarrow {\cal A}$ and 
$\mu' : {\cal A} \otimes {\mbox C} \longrightarrow {\cal A}$ are the 
canonical isomorphisms, defined by 
$$\mu(k \otimes u) = ku = \mu'(u \otimes k), \qquad \forall u \in {\cal A}, 
  \quad \forall k \in {\mbox C} $$

The algebra ${\cal A}$ with the coproduct and the counit has a structure 
of bi-algebra. One extends the algebra ${\cal A}$ by including inverse of 
$x$ which obeys 
$$ x x^{-1} = 1 =  x^{-1} x.$$
If we extend the algebra ${\cal A}$ by adding the inverse of $x$ 
then the algebra ${\cal A}$ admits a C-algebra antihomomorphism 
(coinverse) 
$S: {\cal A} \longrightarrow {\cal A}$ defined by 
$$S(x) = x^{-1}, \qquad S(y) = - x^{-1} y. \eqno(6)$$
Although the coinverse has some properties of an inverse, we have 
$S^2 \neq 1$. Indeed, 
$$S^{-1}(x) = S(x), \qquad S^{-1}(y) = q^{-1}S(y). \eqno(7)$$
The coinverse $S$ satisfies 
$$m \circ (S \otimes \mbox{id}) \circ \Delta = \epsilon = 
  m \circ (\mbox{id} \otimes S) \circ \Delta \eqno(8)$$
where $m$ stands for the algebra product 
${\cal A} \otimes {\cal A} \longrightarrow {\cal A}$. 

We shall demand that 
$$\Delta(1) = 1 \otimes 1. \eqno(9)$$
Then 
$$\epsilon(1) = 1, \qquad S(1) = 1. \eqno(10)$$
Note that 
$$\Delta(x^{-1}) = x^{-1} \otimes x^{-1}. \eqno(11)$$

The coproduct, counit and coinverse which are specified above supply 
the algebra ${\cal A}$ with a Hopf algebra structure. 

\noindent
{\bf 3. Differential Hopf algebra} 

\noindent
{\bf 3.1 Differential algebra} 

\noindent
In this section, we shall give a Hopf algebra structure of the differential 
calculus on the $q$-plane. 
We first note that the properties of the 
exterior differential {\sf d}. The exterior differential {\sf d} is an 
operator which gives the mapping from the generators of ${\cal A}$ to 
the differentials: 
$${\sf d} : u \longrightarrow {\sf d}u, \qquad u \in \{x,y\}. \eqno(12)$$
We demand that the exterior differential {\sf d} has to satisfy two 
properties: the nilpotency 
$${\sf d}^2 = 0 \eqno(13)$$
and the Leibniz rule 
$${\sf d}(f g) = ({\sf d} f) g + f ({\sf d} g). \eqno(14)$$

Let us remind the deformed differential calculus on the quantum plane $^7$ 
$$ \Gamma_1 : \qquad x~ {\sf d}x = q^{-1} {\sf d}x~ x, \qquad 
   x~ {\sf d}y = {\sf d}y~ x, \eqno(15)$$
$$ y~ {\sf d}x = q^{-1} {\sf d}x~ y + (q^{-1} - 1) {\sf d}y~ x, 
   \qquad y~ {\sf d}y = q^{-1} {\sf d}y~ y. $$

The commutation relations between the differentials have the form 
$$ \Gamma_2 : \qquad ({\sf d}x)^2 = 0, \qquad 
  {\sf d}x {\sf d}y = - {\sf d}y {\sf d}x, \qquad 
  ({\sf d}y)^2 = 0.  \eqno(16)$$

A differential algebra on an associative algebra ${\cal B}$ is a 
graded associative algebra $\Gamma$ equipped with an operator {\sf d} 
that has the properties (12)-(14). Furthermore, the algebra $\Gamma$ has to be 
generated by $\Gamma^0 \cup \Gamma^1 \cup \Gamma^2$, where $\Gamma^0$ is 
isomorphic to ${\cal B}$. For ${\cal B}$ we write ${\cal A}$. 
Let $\Gamma$ be the quoitent algebra of the free associative algebra on 
the set $\{x,y,{\sf d}x,{\sf d}y\}$ modulo the ideal $J$ that 
is generated by the relations (1), (15) and (16). 

\noindent
{\bf 3.2 Hopf algebra structure on $\Gamma$}

\noindent
We first note that consistency of a differential calculus with commutation 
relations (1) means that the algebra $\Gamma$ is a graded associative algebra 
generated by the elements of the set $\{x,y,{\sf d}x,{\sf d}y\}$. 
So, it is sufficient only describe the actions of co-maps on the subset 
$\{{\sf d}x,{\sf d}y\}$. 

We consider a map 
$\phi_R : \Gamma \longrightarrow \Gamma \otimes {\cal A}$ such that 
$$\phi_R \circ {\sf d} = ({\sf d} \otimes \mbox{id}) \circ \Delta. \eqno(17)$$
Thus we have 
$$\phi_R({\sf d}x) = {\sf d}x \otimes x $$
$$\phi_R({\sf d}y) = {\sf d}y \otimes 1 + {\sf d}x \otimes y. \eqno(18)$$
We now define a map $\Delta_R$ as follows: 
$$\Delta_R(u_1 {\sf d}v_1 + {\sf d}v_2 u_2) = 
  \Delta(u_1) \phi_R({\sf d}v_1) + \phi_R({\sf d}v_2) 
  \Delta(u_2). \eqno(19)$$
Then it can be checked that the map $\Delta_R$ leaves invariant 
the relations (15) and (16). One can also check that the following 
identities are satisfied: 
$$(\Delta_R \otimes \mbox{id}) \circ \Delta_R = 
  (\mbox{id} \otimes \Delta) \circ \Delta_R, 
  \qquad (\mbox{id} \otimes \epsilon) \circ \Delta_R = \mbox{id}. \eqno(20)$$
But we do not have a coproduct for the differential algebra because the map 
$\phi_R$ does not gives an analog for the derivation property (14), 
yet. So we consider another map 
$\phi_L : \Gamma \longrightarrow {\cal A} \otimes \Gamma$ such that 
$$\phi_L \circ {\sf d} = (\mbox{id} \otimes {\sf d}) \circ \Delta \eqno(21)$$
and a map $\Delta_L$ with again (19) by replacing $L$ with $R$. 
The map $\Delta_L$ also leaves invariant the relations (15) and (16), 
and the following identities are satisfied: 
$$(\mbox{id} \otimes \Delta_L) \circ \Delta_L = 
  (\Delta \otimes \mbox{id}) \circ \Delta_L, 
  \qquad (\epsilon \otimes \mbox{id}) \circ \Delta_L = \mbox{id}. \eqno(22)$$

To denote the coproduct, counit and coinverse which will be defined on the 
algebra $\Gamma$ with those of ${\cal A}$ may be inadvisable. For this reason, 
we shall denote them with a different notation. Let us define the 
map $\hat{\Delta}$ as 
$$\hat{\Delta} = \Delta_R + \Delta_L \eqno(23)$$
which will allow us to define the coproduct of the differential algebra. 
We denote the restriction of $\hat{\Delta}$ to the algebra ${\cal A}$ by 
$\Delta$ and the extension of $\Delta$ to the differential algebra $\Gamma$ 
by $\hat{\Delta}$. It is possible to interpret the relation 
$$\hat{\Delta}\vert_{\cal A} = \Delta \eqno(24)$$
as the definition of $\hat{\Delta}$ on the generators of ${\cal A}$ and 
(23) as the definition of $\hat{\Delta}$ on differentials. 
One can see that $\hat{\Delta}$ is a coproduct for the differential 
algebra $\Gamma$. 

It is not difficult to verify the following conditions: 

{\bf a)} $\Gamma$ is an ${\cal A}$-bimodule, 

{\bf b)} $\Gamma$ is an ${\cal A}$-bicomodule with left and right coactions 
$\Delta_R$ and $\Delta_L$, respectively, making $\Gamma$ a left 
and right ${\cal A}$-comodule with (20) and (22), and 
$$(\Delta_L \otimes \mbox{id}) \circ \Delta_R = 
  (\mbox{id} \otimes \Delta_R) \circ \Delta_L \eqno(25)$$
which is the ${\cal A}$-bimodule property. So, the triple 
$(\Gamma, \Delta_L, \Delta_R)$ is a bicovariant bimodule over Hopf algebra 
${\cal A}$. In additional, since 

{\bf c)} $(\Gamma, {\sf d})$ is a first order differential calculus over 
${\cal A}$, and 

{\bf d)} {\sf d} is both a left and a right comodule map, i.e. for all 
$u \in {\cal A}$ 
$$(\mbox{id} \otimes {\sf d}) \circ \Delta(u) = \Delta_L({\sf d} u), \qquad 
  ({\sf d} \otimes \mbox{id}) \circ \Delta(u) = \Delta_R({\sf d} u), 
  \eqno(26)$$
the quadruple 
$(\Gamma, d, \Delta_L, \Delta_R)$ is a first order bicovariant differential 
calculus over Hopf algebra ${\cal A}$ as noted in ref. 6. 
 
Now let us return Hopf algebra structure of $\Gamma$. If we define a counit 
$\hat{\epsilon}$ for the differential algebra as 
$$\hat{\epsilon} \circ {\sf d} = {\sf d} \circ \epsilon = 0 \eqno(27)$$
and 
$$\hat{\epsilon}\vert_{\cal A} = \epsilon, \qquad 
  \epsilon\vert_\Gamma = \hat{\epsilon}. \eqno(28)$$
we have 
$$\hat{\epsilon}({\sf d}x) = 0, \qquad \hat{\epsilon}({\sf d}y) = 0 
  \eqno(29)$$
where 
$$\hat{\epsilon}(u_1 {\sf d}v_1 + {\sf d}v_2 u_2) = 
  \epsilon(u_1) \hat{\epsilon}({\sf d}v_1) + 
  \hat{\epsilon}({\sf d}v_2) \epsilon(u_2). \eqno(30)$$
Here we used the fact that ${\sf d}(1) = 0$. 

As the next step we obtain a coinverse $\hat{S}$. For this, it suffices to 
define $\hat{S}$ such that 
$$\hat{S} \circ {\sf d} = {\sf d} \circ S \eqno(31)$$
and 
$$\hat{S}\vert_{\cal A} = S, \qquad 
  S\vert_\Gamma = \hat{S} \eqno(32)$$
where 
$$\hat{S}(u_1 {\sf d}v_1 + {\sf d}v_2 u_2) = 
  \hat{S}({\sf d}v_1) S(u_1) + S(u_2) \hat{S}({\sf d}v_2). \eqno(33)$$
Thus the action of $\hat{S}$ on the generators ${\sf d}x$ and 
${\sf d} y$ is as follows: 
$$\hat{S}({\sf d}x) = - x^{-1} ~{\sf d}x ~x^{-1}, $$ 
$$\hat{S}({\sf d} y) = x^{-1} ~{\sf d}x ~x^{-1} y - x^{-1} ~{\sf d}y. 
  \eqno(34)$$
Note that it is easy to check that $\hat{\epsilon}$ and $\hat{S}$ leave 
invariant the relations (15) and (16). 
Consequently, we can say that the structure 
$(\Gamma, \hat{\Delta}, \hat{\epsilon}, \hat{S})$ is a graded Hopf algebra. 

\noindent
{\bf 4. Hopf algebra structure of forms on ${\cal A}$}

\noindent
In this section we shall define two forms using the generators of 
$\cal A$ and show that the algebra of forms is a graded Hopf algebra. 

If we call them $\theta$ and $\varphi$ then one can define them as follows: 
$$\theta = {\sf d}x ~x^{-1}, \qquad 
  \varphi = {\sf d}y - {\sf d}x ~x^{-1} y. \eqno(35)$$
We denote the algebra of forms generated by two elements $\theta$ and 
$\varphi$ by $\Omega$. The generators of the algebra $\Omega$ with the 
generators of $\cal A$ satisfy the following rules 
$$x \theta = q^{-1} \theta x, \qquad 
  y \theta = q^{-1} \theta y + (q^{-1} - 1) \varphi, $$
$$x \varphi = \varphi x, \qquad y \varphi = \varphi y. \eqno(36)$$
The commutation rules of the generators of $\Omega$ are 
$$\theta^2 = 0, \qquad \theta \varphi = - q \varphi \theta, 
  \qquad \varphi^2 = 0. \eqno(37)$$

We make the algebra $\Omega$ into a graded Hopf algebra with the following 
co-structures: the coproduct 
$\Delta: \Omega \longrightarrow \Omega \otimes \Omega$ is defined by 
$$\Delta(\theta) = \theta \otimes 1 + 1 \otimes \theta, $$
$$\Delta(\varphi) = \varphi \otimes 1 + x \otimes \varphi - y \otimes \theta. 
  \eqno(38) $$
The counit $\epsilon: \Omega \longrightarrow {\mbox C}$ is given by 
$$\epsilon(\theta) = 0, \qquad \epsilon(\varphi) = 0 \eqno(39)$$
and the coinverse $S: \Omega \longrightarrow \Omega$ is defined by 
$$S(\theta) = - \theta, \qquad 
  S(\varphi) = - q^{-1} \varphi x^{-1} - \theta x^{-1} y. \eqno(40)$$
One can easy to check that (3), (5) and (8) are satisfied. Note that 
the commutation relations (36) and (37) are compatible with $\Delta$, 
$\epsilon$ and $S$, in the sense that 
$\Delta(x \theta) = q^{-1} \Delta(\theta x)$, and so on.

\noindent
{\bf 6. Discussion}

\noindent
In this paper, we have constructed a Hopf algebra structure on an algebra 
of differential forms on the extended $q$-plane. 

We have choose the coordinate $x$ of the $q$-plane to be group-like: 
$$\Delta(x) = x \otimes x, \quad S(x) = x^{-1} \eqno(41)$$
whereas for $y$ we have defined 
$$\Delta(y) = y \otimes 1 + x \otimes y, \quad S(y) = - x^{-1} y. \eqno(42)$$
However, we know that 
$$\Delta(y) = y \otimes 1 + 1 \otimes y, \quad S(y) = - y \eqno(43)$$
for the elements $y \in g$ of $U(g)$ for a Lie algebra $g$. So the definition 
(42) is a twisting with $x$ of (43). For this reasons the geometry 
which is studied in this paper is actually more similar to the geometry 
of a twisted cylinder that to the $q$-plane. 

Furthermore, if the coordinates of the extended $q$-plane is defined as 
$$x = q^{H/2}, \quad y = X \eqno(44)$$
then we have 
$$[H, X] = 2 X, \eqno(45)$$
where $[,]$ denotes the Lie bracket. In this basis, the Hopf algebra 
structure takes the following form 
$$\Delta(H) = H \otimes {\bf 1} + {\bf 1} \otimes H, \quad S(H) = - H, \quad 
  \epsilon(H) = 0, $$
$$\Delta(X) = X \otimes {\bf 1} + q^{H/2} \otimes X, \quad 
  S(X) = - q^{- H/2} X, \quad \epsilon(X) = 0. \eqno(46)$$
In that case the algebra ${\cal A}$ can be identified with the subalgebra 
of $U_q(sl(2))$. 

\noindent
{\bf Acknowledgement}

\noindent
This work was supported in part by {\bf T. B. T. A. K.} the 
Turkish Scientific and Technical Research Council. 

We would like to express our deep gratitude to the referee for 
clarifying of the fact in Discussion and for general suggestions on the 
manuscript. 


\end{document}